\documentclass[journal]{IEEEtran}

\usepackage{graphicx}
\usepackage{amssymb, amsmath}

\newtheorem{theorem}{Theorem}[section]

\newtheorem{lemma}[theorem]{Lemma}
\newtheorem{Remark}[theorem]{Remark}
\newtheorem{Proposition}[theorem]{Proposition}

\begin{document}
%
\title{Codes Defined By Forms Of Degree 2\\ On Quadric Surfaces}
\author{Fr\'ed\'eric~A. B. Edoukou
\thanks{Manuscript received November 25, 2005; revised }
\thanks{Fr\'ed\'eric A. B. Edoukou is with the Institut de Math\'ematiques de Luminy, CNRS, Luminy case 907  - 13288  Marseille Cedex 9 - France  E.mail : edoukou@iml.univ-mrs.fr}}

\maketitle

\IEEEpeerreviewmaketitle



%

\begin{abstract}
We study the functional codes $C_2(X)$ defined on projective varieties $X$, in the case where $X\subset \mathbb{P}^3$ is a 1-degenerate quadric or a non-degenerate quadric (hyperbolic or elliptic).  We find the minimum distance of these codes, the second weight, and the third weight.  We also show the geometrical structure of the first weight and second weight codewords. One result states that the codes $C_2(X)$ defined on the elliptic quadrics are good codes according to the table of  A. E. Brouwer.
\end{abstract}
\section{Introduction}
The code $C_2(X)$ over $\mathrm{F}_{q}$ where $X$ is a quadric of $\mathbb{P}^{n}(\mathbb{F}_q)$ the projective space of dimension $n$ has been studied by Y. Aubry $\lbrack 1 \rbrack$. In order to find the minimum distance of this code, Aubry determined some bounds for $N$, the number of intersection points of two quadrics $\mathcal{Q}$ and $\mathcal{Q}^{\prime}$ which are exactly known in the case where one of them is an hyperplane, the case both are proportionnal, and the case where the two quadrics have a common factor of degree 1.\\
When the two quadrics have no common linear factor, he proved that: \\ $N \le \pi_{n-2}+7\frac{q^{n-1}}{q-1}-6\frac{q^{n-2}}{q-1} \quad (q\ge 7)$. Thus, even in the case of  $PG(3,q)$, when the quadric surfaces $X$ are 1-degenerate (cones), or non-degenerate (elliptic or hyperbolic) Aubry was not able to determine precisely the minimum distance of the code $C_2(X)$. His result leads to an upper bound: \\$N \le 2(4q+1)+\frac{1}{q-1}$ which is not the best possible. From the result of Lachaud $\lbrack 2 \rbrack$ proposition 2.3, we deduce that this number of intersection points is bounded above by: $N\le 4(q+1)$.\\
Aubry's result has been improved by D. B. Leep and L. M. Schueller $\lbrack 3 \rbrack$ corollary 5.4 p.172 in the case, where the pair of quadrics cannot be written in a fewer number of variables: if $n+1\ge 4$ and even, then $$N \le  (2q^{n-1}-q^{n-2}+2q^{\frac{n+1}{2}}-3q^{\frac{n-1}{2}}+q^{\frac{n-3}{2}}-1)/(q-1).$$ Despite this improvement, the problem of the minimum distance of the code $C_2(X)$ in  $PG(3,q)$ is only  solved for the hyperbolic quadric.\\
In this paper we find some bounds for the number of zeros of a pair of quadrics which are better than the ones known and we give the exact value of the minimum distance of the code $C_2(X)$ with $X$ respectively a 1-degenerate quadric, or an elliptic quadric. For each one of the three types of quadrics mentioned above, we will determine the geometrical structure of the codewords reaching the first weight; we will compute the second and the third weight, and show the geometrical structure of codewords reaching the second weight. Let us note that the codes defined on the elliptic quadrics are good codes according to the table of Andries E. Brouwer $\lbrack 4 \rbrack$.
\section{Notations and terminology}
We denote by $\mathbb{F}_q$ the field with $q$ elements $(q\ge 3)$. We use the term forms of degree two to describe homogeneous polynomials $f$ of degree two, and $\mathcal{Q}=Z(f)$ the zeros of $f$ in the projective space $\mathbb{P}^{n}( \mathbb{F}_q)$, is a quadric.\\
Let $\mathcal{Q}^{\prime}$ be a quadric, the rank of $\mathcal{Q}^{\prime}$ denoted $r(\mathcal{Q}^{\prime})$, is the smallest number of indeterminates appearing in $f^{\prime}$ under any change of coordinate system. The quadric $\mathcal{Q}^{\prime}$ is said to be degenerate if $r(\mathcal{Q}^{\prime})<n+1$; otherwise it is non-degenerate.
The quadrics of $PG(3,q)$ fall into six orbits  under $PGL(4,q)$ (see lemma 15.3.1 $\lbrack 5 \rbrack$). 
In $\lbrack 5 \rbrack$ p.14 we have Table I of the quadrics in $PG(3,q)$. \\
\begin{table}
\renewcommand{\arraystretch}{1.3}
\caption{Quadrics in $\mathbb{P}^{3}(\mathbb{F}_q)$}
\centering
\begin{tabular}{|c|c|c|}

	\hline
	Rank & Description & $\vert \mathcal{Q} \vert$  \\
  	\hline
	\hline
	1 		      & repeated plane 			  &$ q^2+q+1 $				 							       						  \\
	\hline
	2  & pair of distinct planes	  	  & $ 2q^2+q+1 $					 							       					\\
	
	\hline		
	2 		      & line 			  & $ q+1 $				 							       		 \\
	\hline
	3	      & quadric cone 			  & $ q^2+q+1$				 							       			\\
	\hline
	 	4	       & hyperbolic quadric 			  & $ (q+1)^2$				 							       				 \\
	\hline

	 4             &elliptic quadric                               &$q^2+1$                                                         
	
	                                    \\
 \hline            
\end{tabular}
\end{table}
A regulus is the set of transversals of three skew lines in $\mathbb{P}^3$. It consists of $q+1$ skew lines. Thus, a hyperbolic quadric $\mathcal{Q}$ is a pair of complementary reguli  $(\mathcal{Q}=\mathcal{R}_1\cup \mathcal{R}_2$).\\
 Let $F$ be the vector space of forms of degree 2, and $\vert X\vert$ the number of rational points of $X$ over $\mathbb{F}_q$. We denote by $W_i$ the set of points with homogeneous coordinates $(x_0:...:x_3)\in \mathbb{P}^3$ such that $x_j=0$ for $j<i$ and $x_i=1$. The family $\{W_i\}_{0\le i\le 3}$ is a partition of $\mathbb{P}^3$.
The code $C_2(X)$ is the image of the linear map  
 $c: F
  \longrightarrow
  \mathbb{F}_{q}^{\vert X\vert}$, defined by $c(\mathcal{Q}^{\prime})={(c_x(\mathcal{Q}^{\prime})})_{x\in X}$,  where  $c_x(\mathcal{Q}^{\prime})= \mathcal{Q}^{\prime}(x_0,...,x_3)$ with $x=(x_0,...,x_3) \in W_i$.
 \section{Some general theory}
Now we will recall a general result on non-degenerate quadric surface which can be found in $\lbrack 6 \rbrack$ pp.241-243.\\ Let $\mathcal{Q}$ be a non-degenerate quadric over $\mathbb{F}_{q}$ and $\mathcal{Q}^{\prime}$ an other quadric. We denote by $\overline{\mathcal{Q}}$ the points on $\overline{\mathbb{F}}_{q}$ of the quadric $\mathcal{Q}$ defined over $\mathbb{F}_{q}$. Let us write  $\mathcal{C}=\overline{\mathcal{Q}}  \cap \overline{\mathcal{Q}}^{\prime}$ and $\overline{\mathcal{C}}=\bigcup_{1\le i \le t}\overline{C}_{i}$ the decomposition of $\mathcal{C}$ in irreducible components over $\overline{\mathbb{F}}_q$. Note $d_i$=degree of $\overline{C}_{i}$. Since $\overline{\mathcal{Q}}$ is non-degenerate, then $\overline{\mathcal{Q}}$ is isomorphic to $\mathbb{P}_1\times \mathbb{P}_1$ on $\overline{\mathbb{F}}_{q}$ (see  $ \lbrack 6 \rbrack$  chap.4  example 2 p.241).
The curve  $\overline{C}_{i}$ defines a divisor which we note again $\overline{C}_{i}$ on the non-singular quadric $\overline{\mathcal {Q }}$. The divisor class group of $\overline{\mathcal{Q}}$, denoted by $Cl\overline{\mathcal{Q}}$ is isomorphic to $\mathbb{Z} \oplus \mathbb{Z}$. Note $(m,n)$ the class of $\overline{C}_{i}$. The degree and the arithmetical genus of $\overline{C}_{i}$ are such that:  $$\mathrm{deg }\  \overline{C}_{i}=m+n  \quad(3.1) \quad p_a(\overline{C}_{i})=(m-1)(n-1)\quad(3.2)$$
When $\overline{C}_{i}$ is a line, we get deg $\overline{C}_{i}=1$ and $p_a(\overline{C}_{i})=0$. Except in the case $\overline{C}_{i}$ is a line, any irreducible curve has $m>0$, $n>0$. Using the two relations (3.1) and (3.2) we get an interesting relation between the arithmetical genus of a curve and its degree.  Thus we get the following result. 
\begin{Proposition}
Let $\mathcal{Q} \subset {\mathbb{P}}^3$ be a non-singular quadric surface in ${\mathbb{P}}^3$ and
 $\overline{C}_i$ an irreducible curve on $\overline{\mathcal{Q}}$.\\
\textendash If deg $\overline{C}_{i}=2l+2$ then, $p_a(\overline{C}_{i})\le l^2\quad (\star)$\\
\textendash If deg $\overline{C}_{i}=2l+1$ then, $p_a(\overline{C}_{i})\le l(l-1)\quad (\star \star)$
\end{Proposition}
\begin{lemma}
Let $\overline{l}$ be a line defined over $\overline{\mathbb{F}}_{q}$ in $PG(3, \overline{\mathbb{F}}_{q})$
and $l=\overline{l}\cap PG(3, \mathbb{F}_{q})$. Then $l$ is empty, or a single point or a line defined over ${\mathbb{F}_{q}}$.
\end{lemma}
\begin{proof} 
It is obvious.
\end{proof}
\begin{lemma}
Let $\overline{\mathcal{C}}=\bigcup_{1\le i\le t} \overline{C}_i$
 be the decomposition of a curve $\mathcal{C} \subset \mathbb{P}^{n}(\overline{\mathbb{F}}_{q})$ into irreducible components over $\overline{\mathbb{F}}_q$ and $\sigma \in Gal({\overline{\mathbb{F}}}_q/ \mathbb{F}_q)$. 
Then $\sigma(\overline{C}_i)$ and $\overline{C}_i$ have the same degree. 
If $\overline{C}_i$ is the only component of degree $d_i$, then $\overline{C}_i$ is defined over $\mathbb{F}_q$.
\end{lemma}
\begin{proof}  Let $H$ and $H^{\prime}$ be hyperplanes of $\mathbb{P}^{n}(\overline{\mathbb{F}}_{q})$; deg $\sigma(\overline{C}_{i})= \underset{H^{\prime}\not \supset \sigma(\overline{C}_i)}{\max} \#(\sigma(\overline{C_i})\cap H^{\prime})=  \underset{\sigma^{-1}(H^{\prime})\not \supset \overline{C}_i}{\max} \# \sigma(\overline{C_i}\cap \sigma^{-1}(H^{\prime}))=  \underset{H \not \supset \overline{C}_i}{\max} \#(\sigma(\overline{C}_i\cap H)=   \underset{H \not \supset \overline{C}_i}{\max} \#(\overline{C}_i\cap H)$= deg $\overline{C}_{i}$. 
 The second assertion of the lemma comes from the fact that $\{\overline{C}_{i}\}_{1\le i\le t}$ are conjugate under Gal$({\overline{\mathbb{F}}}_q/ \mathbb{F}_q)$.
 \end{proof}
Let us recall the result of Yves Aubry and Marc Perret on the Weil theorem for singular curves, see $ \lbrack 7\rbrack$. With this result we get a bound for the number of rational points of any irreducible curve not necessary smooth.\\
\begin{Proposition}
Let $C$ be a reduced connected projective algebraic curve over a finite field $\mathbb{F}_{q}$, with $r$ irreducible components and of arithmetical genus $p_a(C)$.
Then   $$ \vert \#C(\mathbb{F}_{q})-(rq+1) \vert \le 2p_a(C)\sqrt{q} $$
\end{Proposition}
In $\lbrack 8 \rbrack$ p.229, we have the following result.
\begin{Proposition}
Let $\overline{C}\subset \mathbb{P}^n(\overline{\mathbb{F}}_{q})$   be an irreducible curve which does not lie in any hyperplane. Then deg $\overline{C}$ $\ge n$.
\end{Proposition}
\section{Intersection of two quadrics in $\mathbb{P}^{3}(\mathbb{F}_q)$}
The quadrics of $PG(2,q)$ fall in four orbits under $PGL(3,q)$ (see Theorem 7.4 $\lbrack 9\rbrack$). In $\lbrack 9\rbrack$, we have the following table of quadrics in $PG(2,q)$.
\begin{table}
\renewcommand{\arraystretch}{1.3}
\caption{Quadrics $\mathcal{Q}$ in $\mathbb{P}^{2}{(\mathbb{F}_q)}$}
\centering
\begin{tabular}{|c|c|c|}

	\hline
	Rank & Description & $\vert \mathcal{Q} \vert$  \\
  	\hline
	\hline
	1 		      & repeated line 			  &$q+1 $				 							       						  \\
	\hline
	2  & pair of distinct lines	  	  & $ 2q+1 $					 							       					\\
	
	\hline		
	2 		      & point			  & $ 1 $				 							       		 \\
	\hline
	3	      &  conic 			  & $q+1$				 							       			\\
	    \hline            
	\end{tabular}
\end{table}
\subsection{Plane Sections Of Quadrics}
We have the following result which can be found in $\lbrack 1\rbrack$ p.7.
\begin{Proposition}
Let $\mathcal{Q}$ be a non-degenerate quadric in $\mathbb{P}^{3}$ and H a plane. Then:\\
(i) If H is tangent  to $\mathcal{Q}$, then $\mathcal{Q}\cap H$ is a degenerate quadric of rank 2 in $\mathbb{P}^2$.\\
(ii) If H is not tangent to $\mathcal{Q}$, then $\mathcal{Q}\cap H$ is a conic in $\mathbb{P}^{2}$.\\
\end{Proposition}
From this proposition, we deduce that for any plane section of a non-degenerate quadric $\mathcal{Q}$ in $\mathbb{P}^{3}$, we get $\vert \mathcal{Q}\cap H\vert\le 2q+1$. We also have the same upper bound for a cone, since three of its $q+1$ lines are not coplanar.
\subsection{Intersection Of A Non-degenerate Quadric With A Pair Of Distinct Planes}
Let $\mathcal{Q}$ be an hyperbolic quadric. The lines fall in four classes: the set of generators, tangents, bisecants, and lines skew to  $\mathcal{Q}$. For $\mathcal{Q}$ elliptic, there is no generator.
The two planes meet at a line denoted by $l$. There are two types of planes with respect to  $\mathcal{Q}$: the tangent planes, and the non-tangent planes (those meeting $\mathcal{Q}$ in a conic). Thus the two planes have the three following configurations: \\
 \subsubsection{each plane is tangent to $\mathcal{Q}$} For $\mathcal{Q}$ hyperbolic, there is no line skew to $\mathcal{Q}$ in each tangent plane; there are only tangent lines, bisecants  or generators. However for a fixed tangent line, there is only one tangent plane passing through it. Thus only bisecants  and generators are intersection of two tangent planes. Therefore we get  $\vert \mathcal{Q}\cap \mathcal{Q}^{\prime} \vert =4q$ or $\vert \mathcal{Q}\cap \mathcal{Q}^{\prime} \vert =3q+1$ according to $l$ is a bisecant or a generator. For $\mathcal{Q}$ elliptic, in each tangent plane there is  no bisecant, there are only tangent lines and lines skew to $\mathcal{Q}$. However for any tangent line, there is only one tangent plane passing through it. Thus only lines skew to $\mathcal{Q}$ are intersection of two tangent planes. Therefore, we get  $\vert \mathcal{Q}\cap \mathcal{Q}^{\prime} \vert =2$. \\
 \subsubsection{one of the plane is tangent and the second is non-tangent to $\mathcal{Q}$} For $\mathcal{Q}$ hyperbolic, $l$ is either a tangent or a bisecant; and this gives respectively $\vert \mathcal{Q}\cap \mathcal{Q}^{\prime} \vert= 3q+1$, or $\vert \mathcal{Q}\cap \mathcal{Q}^{\prime} \vert= 3q$. 
For $\mathcal{Q}$ elliptic, $l$ is either a tangent line or skew to $\mathcal{Q}$; thus we get $\vert \mathcal{Q}\cap \mathcal{Q}^{\prime} \vert \le q+2$.\\
 \subsubsection{each one of the two planes is non-tangent to $\mathcal{Q}$} In this case the line $l$ is not a generator;  $l$ is tangent, bisecant or skew to $\mathcal{Q}$, this give respectively $\vert \mathcal{Q}\cap \mathcal{Q}^{\prime} \vert =2q+1$, $\vert \mathcal{Q}\cap \mathcal{Q}^{\prime} \vert =2q$ or $\vert \mathcal{Q}\cap \mathcal{Q}^{\prime} \vert =2(q+1)$.
\subsection{Intersection Of A Cone With A Pair Of Distinct Planes}
Let $\mathcal{Q}$ be a quadric cone. We know from $\lbrack 5 \rbrack$ p.14 that there are four types of planes with respect to $\mathcal{Q}$:  
the set of tangent planes, planes through a pair of generators, planes through the vertex $S$ containing no generator, planes meeting $\mathcal{Q}$ in a conic.
The lines  fall in six classes:\\
$ \mathcal{L}=\mathcal{L}_1\cup \mathcal{L}_1^{+}\cup \mathcal{L}_1^{-}\cup \mathcal{L}_2^{+} \cup \mathcal{L}_2^{-}\cup \mathcal{L}_3$ when $q$ is  odd. \\
$\mathcal{L}=\mathcal{L}_1\cup \mathcal{L}_1^{\prime}\cup \mathcal{L}_1^{\prime \prime}\cup \mathcal{L}_2^{+} \cup \mathcal{L}_2^{-}\cup \mathcal{L}_3$ when $q$ is even\\ and where $\mathcal{L}_3$, $\mathcal{L}_{1}$, $\mathcal{L}_{1}^{+}$, $\mathcal{L}_1^{-}$, $\mathcal{L}_2^{+}$, $\mathcal{L}_2^{-}$, $\mathcal{L}_{1}^{\prime}$ and $\mathcal{L}_{1}^{\prime \prime}$ are respectively the set of generators, the set tangents at the simple points (excluding the generators),  external vertex tangents which are intersection of pairs of tangent planes, internal vertex tangents which are the remaining lines off $\mathcal{Q}$ through $S$, the set of bisecants of $\mathcal{Q}$,  the set of lines skew to $\mathcal{Q}$, the nuclear line, and the set of lines through $S$ which are neither generators nor the nuclear line. \\
Let  $\mathcal{Q}^{\prime}=H_1\cup H_2$, where $H_1$, $H_2$ are two distinct planes and  $l=H_1\cap H_2$. 
Thus, $\vert \mathcal{Q}\cap \mathcal{Q}^{\prime} \vert$ takes the following values:\\
 \textbf{(i)} $\vert \mathcal{Q}\cap \mathcal{Q}^{\prime} \vert =4q+1$, if each plane contains a pair of generators with $l\in \mathcal{L}_1^{\prime \prime}$ for q even, or $l\in \mathcal{L}_1^{+}$, $l\in \mathcal{L}_1^{-}$  for q odd. For example the cone defined by $f_1= (x_0+x_1)x_2+x_2^2$ and $f_2=x_0x_1$ have exactly $4q+1$ common points.\\
 \textbf {(ii)} $\vert \mathcal{Q}\cap \mathcal{Q}^{\prime} \vert =3q+1$, if each plane contains a pair of generators with $l\in \mathcal{L}_3$, or one of the two planes is tangent to $\mathcal{Q}$ and the second plane contains a pair of generators with $l\in \mathcal{L}_1^+$ or $l\in \mathcal{L}_1^{\prime \prime}$\\
 \textbf{(iii)} $\vert \mathcal{Q}\cap \mathcal{Q}^{\prime} \vert \le 3q$, otherwise. 
\subsection{Some Line Sections Of Quadrics}
 \subsubsection{ Section of Quadrics By A Quadric $\mathcal{Q}$ Of Rank 2 And $\mathcal{Q}$ Is A Line} It is obvious that $\vert \mathcal{Q}\capÊ\mathcal{Q}^{\prime}\vert \le \vert \mathcal{Q} \vert$, therefore $\vert \mathcal{Q}\capÊ\mathcal{Q}^{\prime}\vert \le q+1$.\\
 \subsubsection{ Intersection Of An Elliptic Quadric With A Cone Or A Hyperbolic Quadric}  Since there is no line in the elliptic quadric  $\mathcal{Q}$,  and each  of the $q+1$ lines of the cone, or of one regulus meets $\mathcal{Q}$ in at most  two points, we get $\vert \mathcal{Q}\capÊ\mathcal{Q}^{\prime}\vert \le 2(q+1)$. \\
\subsubsection{ Intersection Of A Cone And A Hyperbolic Quadric}
We will distinguish two cases according to this intersection contains a line or not.\\
\textendash If $\mathcal{Q}$ and $\mathcal{Q}^{\prime}$ have no common line:
 we get $\vert \mathcal{Q}\cap \mathcal{Q}^{\prime} \vert \le 2(q+1)$.\\
\textendash if $\mathcal{Q}$ and $\mathcal{Q}^{\prime}$ have at least one common line:
from the study of the hyperbolic quadric in $\lbrack 5 \rbrack$ pp.23-26, we have the following result.
\begin{Proposition}
If C is any point in $\mathcal{Q}^{\prime}$, there pass exactly two generators (i.e. lines lying in $\mathcal{Q}^{\prime}$) which constitute the intersection with $\mathcal{Q}^{\prime}$ of the tangent plane at $C$. And through $C$, there pass $q+1$ lines lying in the tangent plane at $C$, among them $2$ are generators. The remaining $q-1$ lines through $C$, which lie in the tangent plane meet $\mathcal{Q}^{\prime}$ only in the single point $C$ (there are called tangents to $\mathcal{Q}^{\prime}$). Each of the $q^2$ lines passing through $C$ and not contained in the tangent plane at $C$ intersects $\mathcal{Q}^{\prime}$ in exactly two points, one of which is $C$ (such lines are called bisecants to $\mathcal{Q}^{\prime}$).\\
\end{Proposition}
From this proposition, on the hyperbolic quadric  $\mathcal{Q}^{\prime}$, we can assert first of all that the maximum number of generators the cone $\mathcal{Q }$ can contain is $2$. Secondly, if we suppose that the remaining $q-1$ lines  of the cone passing through the vertex $S$ are all bisecants to  $\mathcal{Q}^{\prime}$, we get 
$ \vert \mathcal {Q } \cap \mathcal{Q}^{\prime} \vert \le3q$.\\
\subsubsection{Intersection of two cones}
Let $\mathcal{C}=\mathcal{Q}\cap\mathcal{Q}^{\prime}$; since the degrees of $\mathcal{Q}$ and $\mathcal{Q}^{\prime}$ are 2 and 2, from the B\'ezout theorem $\lbrack 8\rbrack$ p.227, $\mathcal{C}$ contains at most 4 lines.\\
\textendash if $\mathcal{Q}$ and $\mathcal{Q}^{\prime}$ have no common line, each one of the $q+1$ lines of the cone $\mathcal{Q}^{\prime}$ meeting $\mathcal{Q}$ in at most two points, we get 
$\vert \mathcal{Q}\cap \mathcal{Q}^{\prime} \vert \le 2(q+1)$.\\
\textendash if $\mathcal{Q}$ and $\mathcal{Q}^{\prime}$ have exactly one common line, each one of the remaining $q$ lines of $\mathcal{Q}^{\prime}$ meets $\mathcal{Q}$ in at most two points (the vertex of the cone is one of the two points). Therefore, we get $\vert \mathcal{Q}\cap \mathcal{Q}^{\prime} \vert \le (q+1)+q(2-1)$; thus $\vert \mathcal{Q}\cap \mathcal{Q}^{\prime} \vert \le 2q+1.$\\
\textendash if $\mathcal{Q}$ and $\mathcal{Q}^{\prime}$ have exactly two common lines, each one of the remaining $q-1$ lines of $\mathcal{Q}^{\prime}$ meets $\mathcal{Q}$ in at most two points, we deduce that  $\vert \mathcal{Q}\cap \mathcal{Q}^{\prime} \vert \le 3q$.\\
\textendash if $\mathcal{Q}$ and $\mathcal{Q}^{\prime}$ have at least three common lines: from lemma 3.3, we deduce that  $\mathcal{Q}$ and $\mathcal{Q}^{\prime}$ have  exactly four common lines; thus we get $\vert \mathcal{Q}\cap \mathcal{Q}^{\prime} \vert= 4q+1$. For example, the two cones defined by $f_1=(x_0+x_1)x_2+x_2^2$ and $f_2=(x_2+x_0)x_1+x_1^2$ have exactly four common lines.
\subsection{Intersection Of Two Non-degenerate Quadrics Of The Same Type}
From the result of D. B. Leep and L. M. Schueller $\lbrack 3 \rbrack$ p.172 which gives an interesting upper bound for the number $N$ of intersection points of a pair of quadrics in the affine space, we deduce the following consequence. 
\begin{Proposition}
Let $\{g,h\}$ be a pair of quadratic forms defined over $\mathbb{F}_q$  in $PG(n,q)$ with $g$ non-degenerate.
If $n+1\ge 4$ and even, then: \\$N\le (2q^{n-1}-q^{n-2}+2q^{\frac{n+1}{2}}-3q^{\frac{n-1}{2}}+q^{\frac{n-3}{2}}-1)/(q-1).$\\
\end{Proposition}
From this result, we deduce that in $PG(3,q)$ the intersection of a non-degenerate quadric (hyperbolic, or elliptic) with any quadric is always less than $4q$. 
D. B. Leep and L. M. Schueller have found some pair of quadrics from which the upper bound is attained. Their example includes a case where one of the quadric is hyperbolic. Unfortunately for the elliptic quadric in $PG(3,q)$ they didn't find an example; so the problem of the minimum distance of the code $C_2(X)$ defined on the elliptic quadric remains still unresolved.
In what follows, we will study the section of two non-degenerate quadrics of the same type in order to improve the result of D. B. Leep and L. M. Schueller.\\
\subsubsection{Intersection Of Two Hyperbolic Quadrics}
Let $\mathcal{C}=\overline{\mathcal{Q}}\cap\overline{\mathcal{Q}}^{\prime}$, $\mathcal{C}$ is a curve of degree less than 4 from the B\'ezout theorem $\lbrack 8\rbrack$ p.227. Therefore  $\mathcal{C}$ contains at most 4 lines.\\
\textbf{(i) If $\mathcal{Q}$ and $\mathcal{Q}^{\prime}$ have at most one common line :}  
 Let us consider the regulus $\mathcal{R}_1$  which does not contain the common line.  Each one of the $q+1$ lines of $\mathcal{R}_1$ meets $\mathcal{Q}$ in at most two points. Therefore, we get $\vert \mathcal{Q}\cap \mathcal{Q}^{\prime} \vert \le 2(q+1)$. \\
\textbf{(ii) If $\mathcal{Q}$ and $\mathcal{Q}^{\prime}$ have exactly two common lines}, we get two cases:\\
\textendash the two lines are in the same regulus; since there is no intersection line in the complementary regulus, we deduce that $\vert \mathcal{Q}\cap \mathcal{Q}^{\prime} \vert = 2(q+1)$.\\
\textendash the two lines are in distinct regulus $\mathcal{R}_1$ and $\mathcal{R}_2$; we can consider only the regulus $\mathcal{R}_1$. Each one of the remaining $q$ lines of the regulus $\mathcal{R}_1$ meets $\mathcal{Q}^{\prime}$ in at most two points; we deduce that $\vert \mathcal{Q}\cap \mathcal{Q}^{\prime} \vert \le 3q+1$.\\
 \textbf{(iii) If $\mathcal{Q}$ and $\mathcal{Q}^{\prime}$ have at least three common lines}: From lemma 3.3, we deduce that $\mathcal{Q}$ and $\mathcal{Q}^{\prime}$ have exactly four common lines. In this way we get the three following configurations: \\
 \textendash the four lines are in the same regulus, \\ 
 \textendash three lines on a regulus and the fourth in the complementary regulus, \\
 \textendash two lines in each regulus.\\
We know that in $PG(3,q)$ the number of lines meeting 3 skew lines is $q+1$; and these $q+1$ lines are also skew. Since the hyperbolic quadric is a pair of complementary reguli, one concludes that if two hyperbolic quadrics have at least 3 skew lines in common, then  they are identical. Therefore the last configuration is for interest. Let us note that the two first configurations can be eliminated from the result of D. B. Leep and L. M. Schueller.
We deduce that $\vert \mathcal{Q}\cap \mathcal{Q}^{\prime} \vert= 4q$. For example the two hyperbolic quadrics defined by the polynomials $f_1=x_0x_1+x_2x_3$ and $f_2=x_3x_0+x_1x_2$ have exactly  four common lines for $q$ odd.\\
\subsubsection{Intersection Of Two Elliptic Quadrics}
Let $\mathcal{C}=\overline{\mathcal{Q}}\cap\overline{\mathcal{Q}}^{\prime}$, $\mathcal{C}$ is a curve of degree less than 4 from the B\'ezout theorem $\lbrack 8\rbrack$ p.227. Let  $\overline{\mathcal{C}}=\bigcup_{1\le i\le t}{\overline{C}}_{i} $ the decomposition of $\mathcal{C}$ in irreducible components, where ${\overline{C}}_{i}$ are defined over $\overline{\mathbb{F}}_{q}$ and ${C}_{i}={\overline{C}}_{i}\cap  \mathbb{P}^3(\mathbb{F}_q)$. Note $d_i$=degree of ${\overline{C}}_{i}$ then  $\sum_{1\le i\le t}d_i=d$. Without loss of generality we can assume that $d=4$. Thus we get the five following configurations:\\\\
 \textbf{(i) $t=1$:} here $\overline{\mathcal{C}}$ is an irreductible curve of degree 4 defined over $\mathbb{F}_q$ from lemma 3.3, and from part ($\star$) of proposition 3.1 we get $p_a(\mathcal{C})\le 1$. Thus we get $\vert \mathcal{Q} \cap \mathcal{Q}^{\prime}\vert \le 1+q+2\sqrt{q}$.\\\\
 \textbf{(ii) $t=2$ with $d=1+3$:} this case is impossible; in fact from lemma 3.3 there is a common line defined over $\mathbb{F}_q$ in the two elliptic quadrics. \\\\
 \textbf{(iii) $t=3$:} here $\overline{\mathcal{C}}= \overline{l}_1\cup \overline{l}_2\cup \overline{C}_3$, where $\overline{l}_1$ and $\overline{l}_2$ are lines and $\overline{C}_3$ is a curve of degree 2 defined over $\overline{\mathbb{F}}_{q}$. From lemma 3.2, since there is no line on the elliptic quadric, we get that $l_1$ and $l_2$ contain each one at most one point over $\mathbb{F}_q$. The curve $\overline{C}_3$ is also defined over $\mathbb{F}_q$  and from part ($\star $) of proposition 3.1, $p_a(\overline{C}_3)= 0$ so that $\vert \overline{C}_3(\mathbb{F}_q)\vert \le 1+q$.  Therefore we get $\vert \mathcal{Q}\cap \mathcal{Q}^{\prime} \vert \le q+3$.\\\\
 \textbf{(iv)} $t=4$: here $\overline{\mathcal{C}}= \overline{l}_1\cup \overline{l}_2\cup \overline{l}_3\cup \overline{l}_4$  with  $\overline{l}_1$, $\overline{l}_2$, $\overline{l}_3$ and $\overline{l}_4$ lines defined over $\overline{\mathbb{F}}_{q}$. From lemma 3.2, we get that $\vert \mathcal{Q}\cap \mathcal{Q}^{\prime} \vert \le 4$.\\\\
 \textbf{(v) $t=2$ with $d=2+2$:} here $\overline{\mathcal{C}}= \overline{C}_1\cup \overline{C}_2$, where $\overline{C}_1$ and $\overline{C}_2$ are curves of degree 2. Let $\sigma \in Gal( \overline{\mathbb{F}}_{q}/\mathbb{F}_{q}$), then $\sigma(\overline{\mathcal{C}})=\overline{\mathcal{C}}$; $\sigma(\overline{C}_1)\subset \overline{\mathcal{C}}$, therefore $\sigma(\overline{C}_1)= \overline{C}_1$ or $\sigma(\overline{C}_1)= \overline{C}_2$; thus we get two cases:\\
\textendash(a) $\sigma(\overline{C}_1)= \overline{C}_1$ and $\sigma(\overline{C}_2)= \overline{C}_2$: in this case the two curves $\overline{C}_1$ and $\overline{C}_2$ are defined over $\mathbb{F}_q$. Therefore for $i=1,2$, we get $\vert C_i(\mathbb{F}_q)\vert \le 1+q+2p_a(C_i)\sqrt{q}$ with $p_a(C_i)=0$ from proposition 3.1. Hence $\vert \mathcal{Q}\cap \mathcal{Q}^{\prime} \vert \le 2(q+1)$.\\
\textendash(b) $\sigma(\overline{C}_1)= \overline{C}_2$ and $\sigma(\overline{C}_2)= \overline{C}_1$: here we get $\sigma^2(\overline{C}_1)=\overline{C}_1$ and $\sigma^2(\overline{C}_2)= \overline{C}_2$; thus we can assert that $\overline{C}_1$ and $\overline{C}_2$ are defined over $\mathbb{F}_{q^2}$. We will use an analogous technique as the one used in $\lbrack 10 \rbrack$ p.221.
 Let us remark that: 
 $${C}_1=\sigma(\overline{C}_1\cap \mathbb{P}^3(\mathbb{F}_q))=\sigma(\overline{C}_1)\cap \sigma(\mathbb{P}^3(\mathbb{F}_q  ))= {C}_2$$ $\quad \mathrm{and}\quad \overline{\mathcal{C}} \cap \mathbb{P}^3(\mathbb{F}_q)={C}_1.$\\\\
From proposition 3.5, $\overline{C}_1$ is a plane curve therefore it is defined by an irreductible form over $\mathbb{F}_{q^2}$. The curve $\overline{C}_1$ is defined by $f(x,y,z)=\sum f_{ijk}x^{i}y^jz^k$ with $f_{ijk}\in \mathbb{F}_{q^2}$. 
Let $\{1,\alpha\}$ be a $\mathbb{F}_{q}$-base of $\mathbb{F}_{q^2}$: $f_{ijk}=r_{ijk}+\alpha s_{ijk}$ where $r_{ijk}$ and $s_{ijk}$ belong in $\mathbb{F}_{q}$. We can now write $f(x,y,z)=R(x,y,z)+\alpha T(x,y,z)$ where $R(x,y,z)=\sum r_{ijk}x^{i}y^jz^k$  and $T(x,y,z)=\sum s_{ijk}x^{i}y^jz^k$.\\
Let us denote by  $Z_{\mathbb{F}_{q}}(G)$ the set of zeros of the polynomial $G$ in the projective space $\mathbb{P}^{2}(\mathbb{F}_{q}$). In this way we get: $C_1=Z_{\mathbb{F}_{q}}(R)\cap Z_{\mathbb{F}_{q}}(T)$ which gives $\vert C_1\vert =\vert Z_{\mathbb{F}_{q}}(R)\cap Z_{\mathbb{F}_{q}}(T)\vert$. Since there is no line in the elliptic quadric, we deduce that there is no line in $Z_{\mathbb{F}_{q}}(R)\cap Z_{\mathbb{F}_{q}}(T)$  and finally  $\vert C_1\vert \le q+1$.\\
Thus, in this case we have $\vert \mathcal{Q}\cap \mathcal{Q}^{\prime} \vert \le q+1$.
\section{Codes $C_2(X)$ defined on the quadric cones }
Here $X$ denotes a cone quadric. The additional hypothesis in $\lbrack 3 \rbrack$ on the order of the pair of quadrics could not be satisfied. We apply the results of  IV. (A, C, D) to  determine the parameters of the functional code $C_2(X)$ defined on the cone quadric. Thus, we get the following results.
\begin{Proposition}
Let $\mathcal{Q}^{\prime}$ be a quadric in $PG(3,q)$ and $X$ the quadric cone. we have: 
$$\# X_{ Z(\mathcal{Q}^{\prime})}(\mathbb{F}_{q} ) = 4q+1, 3q+1 \quad \mathrm{or}\quad \# X_{ Z(\mathcal{Q}^{\prime})}(\mathbb{F}_{q} ) \le 3q$$
\end{Proposition}
\begin{theorem}
The code $C_2(X)$ defined on the quadric cone $X$ is an ${\lbrack n,k,d\rbrack }_{q}$- code with: $$n=q^2+q+1, \quad k=9, \quad d=q(q-3)$$
\end{theorem}
\begin{theorem}
The minimum weight codewords correspond to:\\
\textendash \ quadrics which are union of two planes passing through a pair of generators and the line of intersection of the two planes intersects the cone at a singhe point, \\
\textendash \ quadrics cone containing exactly four lines of $X$.
\end{theorem}
\begin{theorem}
The second and the third weight of the code $C_2(X)$ defined over the quadric cone are respectively $w_2=q(q-2)$, $w_3=(q-1)^2$.
\end{theorem}
\begin{theorem}
The codewords of second weight correspond to:\\
\textendash \ quadrics which are union of two planes, each passing through a pair of generators and the line of intersection of the two planes is a generator,\\  
\textendash \ quadrics which are union of two planes, one tangent to the cone and the second passing through a pair of generators and the line of intersection of the two planes $l\in \mathcal{L}_1^+$ or $l\in \mathcal{L}_1^{\prime \prime}$.
\end{theorem}
\section{Codes $C_2(X)$ defined on non-degenerate quadrics}
Here $X$ denotes a non-degenerate quadric.  
Since the parameters of the code $C_2(X)$ vary according to the type of the quadric $X$, we need to distinguish two cases: hyperbolic quadrics and elliptic quadrics. 
\subsection{The Parameters Of The Codes $C(X)$ Defined On Hyperbolic Quadrics}
As previously indicated, the upper bound in $\lbrack 3 \rbrack $ p.172 is optimal in this case;  it gives the minimum distance. In the other hand it can not describe more than the first weight; this can be done by applying the results of IV. (A-B, D.1-3, E.1). Thus, we have the following results.
\begin{Proposition}
Let $\mathcal{Q}^{\prime}$ be a quadric in $PG(3,q)$ and $X$ a hyperbolic quadric. We have: 
$$\# X_{ Z(\mathcal{Q}^{\prime})}(\mathbb{F}_{q} ) = 4q, 3q+1\quad \mathrm{or}\quad \# X_{ Z(\mathcal{Q}^{\prime})}(\mathbb{F}_{q} ) \le 3q$$
\end{Proposition}
\begin{theorem}
The code $C_2(X)$ defined on the hyperbolic quadric $X$ is an ${\lbrack n,k,d\rbrack }_{q}$- code with: $$n=(q+1)^2, \quad k=9, \quad d=(q-1)^2$$
\end{theorem}
\begin{theorem}
The minimum weight codewords correspond to:\\
\textendash \ quadrics which are union of two tangent planes to $X$ and the line of intersection of the two planes intersects the hyperbolic quadric at two points, \\
\textendash \ hyperbolic quadrics containing exactly four lines of $X$ with two lines  in each regulus of $X$.
\end{theorem}
\begin{theorem}
The second and the third weights of the code $C_2(X)$ defined on the hyperbolic quadric are respectively $w_2=q(q-1)$, $w_3=q^2-q+1$.
\end{theorem}
\begin{theorem}
The codewords of second weight correspond to:\\
\textendash \ quadrics which are union of two tangent planes to $X$ and the line of intersection of the two planes is a generator of $X$, \\
\textendash \ quadrics which are union of a tangent plane, and a non-tangent plane to $X$, and the line of intersection of the two planes is tangent  to $X$,\\
\textendash \ hyperbolic quadrics containing exactly two generators  of $X$ in distinct reguli and the remaining q lines of one regulus are all bisecants to $X$.
\end{theorem}
\subsection{The Parameters Of The Codes $C_2(X)$ Defined On Elliptic Quadrics}
Here the upper bound from the result of D. B. Leep and L. M. Schueller is not the best possible; it  cannot  give the minimum distance. The minimum distance and the firsts weights of the code $C_2(X)$, can exactly be determined by a straight application of IV. (A-B, D.1-2 and E.2). Thus, we get the following results.
\begin{Proposition}
Let $\mathcal{Q}^{\prime}$ be a quadric in $PG(3,q)$ and $X$ an elliptic quadric. We have: 
$$\# X_{ Z(\mathcal{Q}^{\prime})}(\mathbb{F}_{q} ) = 2(q+1), 2q+1 \quad \mathrm{or}\quad  \# X_{ Z(\mathcal{Q}^{\prime})}(\mathbb{F}_{q} ) \le 2q$$
\end{Proposition}
\begin{theorem}
The code $C_2(X)$ defined on the elliptic quadric $X$ is an ${\lbrack n,k,d\rbrack }_{q}$- code with: $$n=q^2+1, \quad k=9, \quad d=q^2-2q-1.$$
\end{theorem}
\begin{Remark}
 These functional codes $C_2(X)$ defined on the elliptic quadrics are good codes according to the table of Andries E. Brouwer $\lbrack 4\rbrack$.
\end{Remark}
In fact, for $q=3, 4, 5, 7, 8, 9$ we get respectively the following codes:
${\lbrack 10,9,2\rbrack }_{3}$- code, ${\lbrack 17,9,7\rbrack }_{4}$- code, ${\lbrack 26,9,14\rbrack }_{5}$- code,  ${\lbrack 50,9,34\rbrack }_{7}$- code, ${\lbrack 65,9,47\rbrack }_{8}$- code, ${\lbrack 82,9,62\rbrack }_{9}$- code.
\begin{theorem}
The minimum codewords correspond to:\\
\textendash \ quadrics which are union of two non-tangent planes to $X$ and the line of intersection of the two planes is skew to $X$, \\
\textendash \ quadrics cones which vertex doesn't lie on $X$, and all its $q+1$ lines are bisecants to $X$,\\
\textendash \ hyperbolic quadrics with all the $q+1$ lines of a regulus are bisecants to $X$.
\end{theorem}
\begin{theorem}
The second and third weights of the code $C_2(X)$defined on the elliptic quadric $X$ are respectively $w_2=q(q-2)$, $w_3=(q-1)^2$.
\end{theorem}
\section*{Acknowledgment}
The author would like to thank Mr F. Rodier and Mr. D. B. Leep.  Theirs remarks and comments encouraged him to work on the problem. 
\section*{References}
\noindent {\footnotesize \lbrack 1\rbrack   \ Y. Aubry,  Reed-Muller codes associated to projective algebraic varieties. In ``Coding Theory  and Algebraic Geomertry and Coding Theory ". (Luminy, France, June 17-21, 1991).  Lecture Notes in Math., Vol. 1518 pp.4-17, Springer-Verlag, Berlin, 1992.\\  
\lbrack 2\rbrack  \ G. Lachaud, Number of points of plane sections and linear codes defined on algebraic varieties;  in " Arithmetic, Geometry, and Coding Theory " . (Luminy, France, 1993), Walter de Gruyter, Berlin-New York, 1996, pp 77-104. \\ 
\lbrack 3\rbrack \ D. B. Leep and L. M. Schueller, Zeros of a pair of quadric forms defined over finite field. Finite Fields and Their Applications 5,  157-176 (1999).\\ 
\lbrack 4\rbrack \ A. E. Brouwer, Bounds on minimum distance of linear codes,
htpp://www.win.tue.nl/\\\verb+~+aeb/voorlincod.html\\
\lbrack 5\rbrack  \ J. W. P. Hirschfeld, Finite projective spaces of three dimensions, Clarendon press. Oxford 1985. \\ 
\lbrack 6\rbrack \  I. R. Shafarevich, Basic algebraic geometry 1, Springer-Verlag, 1994.\\ 
\lbrack 7\rbrack   \ Y. Aubry and M. Perret, Connected projective algebraic curves over finite fields. Preprint, IML, Marseille, France, 2002.\\ 
\lbrack 8\rbrack  \ J. Harris,  Algebraic Geometry  (A first course), Graduate texts in mathematics 133, Springer-Verlag,\\ 
\lbrack 9\rbrack  \ J. W. P. Hirschfeld, Projective Geometries Over Finite Fields (Second Edition) Clarendon  Press. Oxford 1998.\\ 
\lbrack 10\rbrack \ J. P. Cherdieu and R. Rolland, On the number of points of some hypersurface in $\mathbb{F}^n_q$, Finite Fields and Their Applications 2,  214-224 (1996).}

\end{document}